\documentclass{article}
\usepackage{amsmath}
\usepackage{amsmath,amscd}
\usepackage{amssymb,amsfonts,amsmath,amscd}
\usepackage[margin=0.5in,footskip=0.25in]{geometry}
\usepackage{fancybox}
\usepackage{pgf,pgfarrows,pgfnodes,pgfautomata,pgfheaps,pgfshade}
\usepackage[all,graph,frame]{xy}
\usepackage{appendix}
\usepackage{enumerate}

\def\R{\mathbb{R}}

\def\proof{{\em Proof}}
\def\qed{{$\Box$}}

\newtheorem{theorem}{Theorem}

\newtheorem{lemma}{Lemma}

\newtheorem{proposition}{Proposition}
\newtheorem{remark}{Remark}

\begin{document}


\title{Generalizing Gale's theorem on
backward induction and domination of strategies}
\author{Vladimir Gurvich
\footnote
{National Research University, Higher School of Economics (HSE) Moscow Russia
             \textit{vgurvich@hse.ru}}}
\maketitle

\;\;\;\;\;\;\;\;\;\;\;\;\;  To Yuriy Borisovich Germeyer to his centenary.

\begin{abstract}
In 1953 David Gale noticed that
for every  $n$-person game in extensive form with perfect information modeled by an arborescence
(a rooted tree)  some special Nash equilibrium in pure strategies can be found
by an algorithm of successive elimination of leaves, which is now called
the {\em backward induction}. 
(The result can be easily extended from the trees to the acyclic directed graphs.)
He also noticed the same procedure, performed for the normal form
of this game, turns into successive elimination of dominated strategies of the players
that results in a single strategy profile $(x_1, \ldots, x_n)$,
which is called a {\em domination equilibrium}
(DE)  and appears to be a  Nash-equilibrium (NE) too.
In other words, the game in normal form
obtained from a positional game with perfect information
is {\em dominance-solvable} (DS)  and also {\em Nash-solvable} (NS).
Yet, an arbitrary game in normal form may be not DS.
We strengthen Gale's results as follows.
Consider several successive eliminations of dominated strategies
that begins with  $X = X_1 \times \ldots \times X_n$  and ends in
$X' = X'_1 \times \ldots \times X'_n$.
We will call  $X'$  a D-box of  $X$.
Our  main (although obvious) lemma claims that for any
$i \in I = \{1, \ldots, n\}$  and for any strategy  $x_i \in X_i$
its ``projection" to a D-box  $X'$  of  $X$  is  dominated by a strategy  $x'_i \in X'_i$.
It immediately follows that any DE is an NE and, hence, DS implies NS.
It is enough to apply the lemma
in the case when  $X'$  consists of a single strategy profile.
Also this lemma implies that, in general,
the domination procedure is well-defined in the following sense.
A D-box  $X'$  of  $X$  is called {\em terminal}
if it is domination-free, that is, it contains no pair of strategies
such that one of them is dominated by the other.
Any two  terminal D-boxes $X'$  and  $X''$  of  $X$  are equal.
More precisely, there exist  $n$  permutations
$\pi = (\pi_1, \ldots, \pi_n)$, with  $\pi_i : X_i \rightarrow X_i$
for  $i \in I$, that  transform  $X'$  into  $X''$, that is,
$\pi(X') = X''$  and the payoffs are respected.
\newline
We also recall some published results on dominance-solvable game forms.
\end{abstract}

\section{Main results}
Now, it seems too late to publish these results.
In Spring of 1971, I gave a talk
"Successive elimination of dominated
strategies for games in normal and positional forms"
for the seminar on Game Theory,
guided by my advisor Yuriy Germeyer,
in Computing Center of the USSR Acaddemy of Sciences.
Then, I repeated this presentation many times
on the seminars and lectures in game theory in different places;
yet, it did not appear in writing,
except for \cite[Section 4]{GGZ86}, which is a microfilm,
in Russian, and I strongly doubt that it is still available.

\subsection{Games and game forms; domination is transitive and anti-symmetric}
\label{ss00}
Let  $I = \{1, \ldots, n\}$  be the set of {\em players},
$X_i$  be a finite set of {\em strategies} of $i \in I$,
and  $X = X_1 \times \ldots \times X_n$  be
the direct product of these  $n$  sets.
Then, $x = (x_1, \ldots, x_n) \in X_1 \times \ldots \times X_n  = X$
is called a {\em strategy profile}
and  $u : I \times X \rightarrow \R$  a utility function;
$u(i, x)$ is interpreted as the profit of the player  $i \in I$
in case the profile  $x \in X$  is realized.
The pair  $(X, u)$  defines a game in normal form.

Given two strategies  $x_i, x'_i \in X_i$ of a player
$i \in I$, we say that
$x_i$  dominates  $x'_i$, or equivalently,
$x'_i$  is {\em dominated} by  $x_i$, and write
$x_i \succeq x'_i$  or  $x'_i \preceq x_i$
if  $u(i, x) \geq u(i, x')$
for any two strategy profiles  $x$ and $x'$
that coincide in all coordinates  but  $i$ and the  latter
is  $x_i$  in  $x$   and  $x'_i$  in  $x'$.
Furthermore, we will say that  $x_i$  and  $x'_i$ are equivalent
and write  $x_i \simeq x'_i$
if   $u(i', x) = u(i', x')$  for any  $i' \in I$
(in particular, for  $i' = i$).
By definition, two equivalent strategies dominate each other.
Obviously, domination is transitive and idempotent:

\medskip

$x_i \succeq x'_i$   and  $x'_i \succeq x''_i$  implies
$x_i \succeq x''_i$; \;\;\;
$x_i \succeq x_i$;  \;\;\;  for any  $i \in I$.
\medskip

We would also like it to be anti-symmetric:
$x_i \succeq x'_i$   and  $x'_i \succeq x_i$  implies
$x_i \simeq x'_i$.  
The reason is clear.
We plan to eliminate the dominated strategy  $x'_i$  from  $X_i$, since
it cannot be better than $x_i$.
Yet, it may happen that  $u(i, x) = u(i, x')$,
for all  $x, x'$  considered above, while
$u(j, x) \neq u(j, x')$  for some  $x, x'$  and  $j \in I$.
In other words, the strategies  $x_i$  and  $x'_i$
are equally good for  $i$  in all cases, so
we could eliminate any of them and  $i$  will be indifferent.
Yet, the difference between  $x_i$ and $x'_i$  may be very important for other players.
Thus, without anti-symmetry of domination no uniqueness will hold.
To resolve this issue, let us simply require for all  $x, x' \in X$  that

\medskip

$u(i, x) = u(i, x')$  for all  $i \in I$   whenever
$u(i', x) = u(i', x')$  for some  $i' \in I$,

\medskip
\noindent
if  $x$  and  $x'$  are equally good for a player then
they are equally good for all players,
or in other words, represent the same outcome of the game.

\medskip

Thus, we naturally come to the concepts of a {\em game form},
defined as a mapping  $g : X \rightarrow A$, where
$A  = \{a_1, \ldots, a_p\}$  is a set of outcomes of the game.
Naturally, we assume that two distinct outcomes
$a, a' \in A$  cannot be equally good for all player,
since otherwise, why are they distinct.
Hence, they cannot be equally good for any player, that is,

\medskip

$u(i, a) = u(i, a')$  for  some  $i \in I$  implies that  $a = a'$.

\medskip

Then,  $u = (u_1, \ldots, u_n)$  is a preference  profile
rather than a real valued utility function:
$u_i$  is a complete order over  $A$
reflecting the   preference of player  $i \in I$.
Now, the pair  $(g, u)$  defines a game in normal form;
two strategies  $x_i$  and  $x'_i$  are equivalent,
if  $g(x) = g(x')$  is the same outcome from  $A$
for  $x, x' \in X$  considered above, and
domination becomes anti-symmetric:

\medskip

$x_i \succeq x'_i$  and  $x'_i \succeq x_i$  imply
that  $x_i \simeq x'_i$.

\subsection{Main lemma}
\label{ss01}
Given  $X = X_1 \times \ldots \times X_n$   and
$Y = Y_1 \times \ldots \times Y_n$
we will call  $Y$  a {\em box} of  $X$  and
write  $Y \subseteq X$   if
$Y_i \subseteq X_i$  for all  $i \in I$.
Given a game form $g : X \rightarrow A$, a game  $(g, u)$,
and two strategies  $x_i, x'_i \in X_i$, we say that
$x_i$  dominates  $x'_i$  modulo  $Y$   and write
$x_i \succeq x'_i \mod Y$  if
$u(i, x) \geq u(i, x')$  for every  $x$  and $x'$  such that
$x_j = x'_j  \in  Y_j$  for all  $j \in I \setminus \{i\}$.
To visualize, we can say that
projection of  $x_i$  into  $Y$  dominates
projection of  $x'_i$  into  $Y$.
Obviously, domination respects  projection: if

\medskip

$x_i \succeq x'_i \mod X$  and  $Y \subseteq X$ then
$x_i \succeq x'_i \mod Y$

\medskip

We will need several definitions related to boxes:
\begin{itemize}
\item
A box  $Y \subseteq X$  is  {\em domination-free}
if for any  $i \in I$  it contains no two strategies
$x_i, x'_i \in Y_i$  such that  $x_i \succeq x'_i \mod Y$.
\item
We will successively eliminate dominated strategies
of (different) players  $i \in I$, one by one, in some order.
Any box  $Y \subseteq X$  obtained in several such steps
is called a {\em D-box}.
\item
A domination-free  D-box will be called {\em terminal}.
\item
A terminal D-box that consists of a unique strategy profile
will be called a {\em domination equilibrium} (DE).
\end{itemize}

Note that on each step, the set of dominated strategies
of a player may depend on
which strategy of another player was eliminated on the previous steps,
if these two players are distinct.
Hence, the procedure is not unique, it branches.
Yet, making use of the following simple claim,
we will prove, in Section \ref{ss02} that, in a sense, the resulting
domination-free D-box is well-defined.

\begin{lemma}
\label{l1}
For a D-box  $Y \subseteq X$  and a strategy  $x_i \in X_i$,
there exists a strategy  $y_i \in Y_i$  such that  $y_i \succeq x_i \mod Y$.
\end{lemma}

\proof.
If  $x_i \in Y_i$  then  $x_i \succeq x_i$  and there is nothing to prove.
If  $x_i \not\in Y_i$  then at some step it was dominated
by a strategy  $x'_i \in X_i$.
If  $x'_i \in Y_i$  then just set  $y_i = x'_i$  and we are done.
If  $x'_i \not\in Y_i$  then at some step it was dominated by
a strategy  $x''_i \in X_i$, etc.
Sooner or later we will obtain some  $x^{''\ldots'} \in Y_i$.
Since domination is transitive and it respects projection, our claim follows.
\qed.

From the definition of NE and Lemma \ref{l1}
we immediately derive the  next statement.

\begin{proposition}
\label{p1}
Any domination equilibrium is a Nash equilibrium.
\qed
\end{proposition}

\subsection{All terminal D-boxes are equal}
\label{ss02}
Given a game  $(g, u)$  in normal form.
As we mentioned, typically, one can eliminate
dominated strategies of the players
in many different orders, thus, getting different terminal D-boxes.
Yet, in the following sense, they are all equal.

\begin{theorem}
\label{t1}
For any two terminal D-boxes  $Y, Z \subseteq X$
there exist a profile of  $n$  permutations
$\pi = (\pi_i : X_i \rightarrow X_i \mid i \in I)$  such that
$\pi(X) = Y$  and  $g(x) = g(\pi(x))$  for all  $x \in X$,
\end{theorem}

In particular, Theorem \ref{t1} implies that
all terminal D-boxes are of the same size and
in case it is $1 \times \ldots \times 1$,
all corresponding DE realize the same outcome.

\medskip

\proof \; of Theorem \ref{t1}.
Consider any two sequences of eliminating
dominated strategies in  $X$, one by one,
that result in  $Y$  and  $Z$.
There is a unique step of the second sequence
when a strategy  $x'_i \in Y_i$  is eliminated first time,
by a strategy  $x_i \in X_i$.
Then,  $x_i \succeq x'_i \mod Y$, by projectivity.
Note that  $x_i \not\in Y_i$,  since  $Y$  is a terminal D-box.
By Lemma \ref{l1}, there exists  $x''_i \in Y_i$  such that
$x''_i \succeq x_i \mod Y$.
Then, by transitivity, $x''_i \succeq x'_i \mod Y$.
But  $Y$  is terminal...
Can we stop here just saying that
two strategies of  $Y$  cannot dominate each other,
since  $Y$  is a terminal D-box?
No, because this is true only for two distinct strategies, and
we can only conclude that  $x'_i = x''_i \simeq x_i$.
The equality holds, because  $Y$  is a terminal D-box, while
the equivalence holds, because domination is anti-symmetric.
And ... here we can stop.
Indeed, whenever we eliminate a strategy from  $Y$
an equivalent, $\mod Y$, strategy of the same player is left in the current box.
Replace  $x'_i = x''_i$  by  $x_i$  and proceed.
At the end, when we reach  $Z$,  a required permutation profile,
which transfers  $Y$  to  $Z$  will be obtained.
\qed

\section{Dominance-solvable (and Nash-solvable) game forms}
\label{s1}

\subsection{Backward induction in extensive and normal form}
\label{ss10}
The backward induction (BI) algorithm was developed
for $n$-person games with perfect information
in extensive (positional) form,
modeled by an arborescence (that is, by a rooted tree))
by Gale \cite{Gal53} and  Kuhn \cite{Kuh53};
see also \cite[Theorem 1]{Kuh50}.

For positional games BI works as follows.
Take a position  $v$  every move from which enters a terminal;
chose a move $(v, v')$  maximazing the payoff of the player  $i \in I$
that controls  $v$; assign to  $v$  the whole payoff vector in  $v'$,
make  $v$  a new terminal, and eliminate all terminals succeeding  $v$.
We already mentioned the problem caused by possible ties
and explained how to resolve it.
Repeat the above procedure until the root  $v_0$  becomes the terminal.
Thus, we obtain the resulting payoff vector, and also the strategy profile
$x = (x_1, \ldots, x_n)$, which is a DE.

``It is somewhat surprising, then, that the same solution
could have been obtained by examining
only the normal form of the game" \cite{Gal53}.
Yet, maximizing the payoff of the player  $i \in I$
that controls  $v$  over  all terminals  $a \in A$  that succeed $v$,
we do in fact the following:
In the extensive form, we fix a move  $(v, a)$  maximizing  $u_i$,
assign the whole payoff profile in  $a$   to  $v$, and then,
eliminate all moves  $(v, a)$  making  $v$  a new terminal.
In the normal form, we select a strategy  $x_i$
that maximizes the payoff of  $i$  in  $v$
making a move  $(v, a)$  and some moves in the remaining positions
and a strategy  $x'_i$  that makes the same moves
in the remaining positions and some other move  $(v, a')$  in  $v$.
Obviously, $x_i$  dominates  $x'_i$.
We consider all such comparable pairs, one by one
(typically, there are very many of them)
and each time we eliminate the dominated strategy.

\begin{remark}
\label{r1}
For each such pair  $x_i, x'_i$
either  $g(x) = g(x')$  or
$g(x) = a$  and  $g(x') = a'$, where  $a$  and  $a'$ are
introduced in the previous paragraph,
while for  $x$ and  $x'$  see the definition of domination.
In other words, the strategies  $x_i$  and  $x'_i$  differ only
in outcomes  $a$  and  $a'$,  while all other corresponding outcomes coincide.
\end{remark}

A game form  $g : X \rightarrow A$  is called {\em positional}
if it is generated by an extensive $n$-person game structure
with perfect information modeled by rooted tree.
The previous results immediately imply the next statement.

\begin{proposition}
Positional game forms are dominance-solvable.
\qed
\end{proposition}

Several examples of positional game forms
together with the corresponding positional structures
can be found in \cite{Gur09}  and  \cite{BG15}.
For reader's convenience we reproduce
some positional game forms here.

\medskip

$a_1 a_2 a_3 a_2 a_3 \;\;\;\;\;\; a_1 a_1\;\;\;\;\;\;  a_1 a_1 a_2 a_2 \;\;\;\;\;\;\;\;\;\;\;\;  a_1 a_1 \;\; a_3 a_4$

$a_1 a_4 a_4 a_5 a_5 \;\;\;\;\;\; a_2 a_3\;\;\;\;\;\;  a_3 a_4 a_3 a_4 \;\;\;\;\;\;\;\;\;\;\;\;  a_2 a_2 \;\; a_3 a_4$

$a_1 a_4 a_4 a_6 a_6$

\medskip

The first three are two-person game forms, while the last one is a three-person game form.

\subsection{Characterizing positional game forms}
\label{ss11}
Their characterization was given in \cite[Chapter 5]{Gur78} and \cite{Gur82};
see also \cite{Gur84}, \cite{Gur09}, and \cite{BG15} for more details.
Here we just summarize briefly the results.

Let $g : X = \prod_{i \in I} X_i  \rightarrow A$  be an $n$-person game form.
Standardly, $I$ and $A$  denote respectively the sets
of the players (voters) and the outcomes (candidates) in the game (voting) theory.
We assume that mapping  $g$  is surjective, $g(X) = A$, but
typically it is not injective, in other words,
different strategy profiles may determine the same outcome.
The subsets $K \subseteq I$  and  $B \subseteq A$
are called {\em coalitions} and {\em blocks}, respectively.
A strategy  $x_K = (x_i \; | \; i \in K)$  of a coalition  $K$ is
defined as a collection of strategies, one for each coalitionist of $K$.
A pair of strategies  $x_K$  of coalition $K$ and
$x_{I \setminus K}$ of the complementary coalition $I \setminus K$
uniquely defines a strategy profile  $x = (x_K, x_{I \setminus K})$.

We say that a non-empty coalition  $K$  is {\em effective} for a block  $B$
if this coalition has a strategy such that
the resulting outcome is in  $B$  for any strategies of the opponents;
in other words, if there is a strategy  $x_K$
such that for any strategy  $x_{I \setminus K}$
for the obtained strategy profile  $x = (x_K, x_{I \setminus K}$
we  have  $g(x) \in B$.
We write  $E(K,B) = 1$  if  $K$  is effective for  $B$  and  $E(K,B) = 0$  otherwise.

Since  $g$  is surjective, by the above definition, we have
$E_g(I, B)= 0$  if and only  $B = \emptyset$.
The value  $E_g(K,B)$  was not yet defined for $K = \emptyset$.
By convention, we set  $E_g(\emptyset, B)= 1$  if and only if  $B = A$.

The obtained function   $E_g: 2^I \times 2^A = 2^{I \cup A} \rightarrow \{0,1\}$
is  called the {\em effectivity function} (of the game form $g$).
By definition, $E_g$ is a Boolean function whose set of variables
is the union $I \cup A$ of all players and outcomes.

Note that equalities $E_g(K,B) = 1$ and
$E_g(I \setminus K, A \setminus B)= 1$  cannot hold simultaneously.
Indeed, otherwise there would exist strategies $x_K, x_{I \setminus K}$
and blocks  $B,  A \setminus B$  such that
$g(x) \in B \cap (A \setminus B) = \emptyset$, that is,
$g$  is not defined on the strategy profile  $x = (x_K, x_{I \setminus K})$.
Note that cases  $K = \emptyset$  and  $K = I$  are covered by the above convention.
Thus, we have  $E_g(K,B) = 1 \Rightarrow E_g(I \setminus K, A \setminus B) = 0$.

A game form  $g$  is called {\em tight} if the inverse implication holds too:
$E_g(I \setminus K, A \setminus B)= 0 \Rightarrow E_g(K,B) = 1$, or equivalently, if
$E_g(K,B) = 1 \Leftrightarrow E_g(I \setminus K, A \setminus B)= 0$
for all pairs  $K \subseteq I$  and $B \subseteq A$, or,
in other words, if  $E_g$  is a self-dual Boolean function,
see, e.g.,~\cite[Part 1 Chapter 4]{CH11} for the definition.
A game form  $g$  is called {\em weakly tight}
if the required implication holds when
$K$  or  $I \setminus K$  is a single player,
that is, $|K| = 1$  or   $|K| = n - 1$.

A game form  $g$  is  called {{\em rectangular} if
$g(x') = g(x'') = a$  implies that  $g(x) = a$  for any
$x = (x_i \mid i \in I)$  such that
$x_i = x'_i$  or  $x_i = x'_i$   for all  $i in I$, where
$x' = (x'_i \mid i \in I)$  and $x'' = (x''_i \mid i \in I)$;
in other words, mapping  $g : X \rightarrow A$
assigns a box in $X$  to each outcomes  $a \in A$.

\begin{theorem} (\cite[Chapter 5]{Gur78}).
The following three property of a  game form  $g$   are equivalent:

$g$  is positional,
$g$  is tight and rectangular,
$g$  is weakly tight and rectangular.
\end{theorem}

The result was announced in \cite[Remark 3]{Gur75}; see also \cite{Gur82}, \cite{Gur09}, \cite{GG11} and \cite{BG15}).

We leave to the reader verify these three properties for the four game forms  given above.

\subsection{Tight, totally tight, Nash-solvable, and acyclic game forms for $n = 2$  and  $n > 2$}
\label{ss12}
For the two-person game forms, $n = 2$, the following chain of important implications was shown in \cite{BGMP10}.
Moreover, it was also shown that no other implication may hold between these five properties.

\medskip

$TT \Leftrightarrow AC \Rightarrow DS \Rightarrow NS \Leftrightarrow T$

\begin{itemize}
\item {\em Acyclicity} (AC) means the absence of improvement cycles for any  $u$.
\item {\em Total tightness} (TT) means that
every $2  \times 2$  box contains a constant (that is, containing only one outcome) line.
\end{itemize}

We refer the reader to   \cite{BGMP10}  and  \cite{BCG12} for more details
and,  in particular, for extending the concept of  TT  for  $n > 2$.
The last three properties were already defined above for arbitrary  $n$.

Tightness and NS are equivalent for the two-person game forms.
This result was obtained in \cite{Gur75} and repeated in a more
focused paper \cite{Gur88}; see also
\cite[Chapter 2]{Gur78} and much more recent
\cite{LRMDS17}, \cite{Gur17}, \cite{GK17}.

Yet, already for $n = 3$
tightness is not necessary and not sufficient for NS;
see examples in \cite[Remark 3]{Gur75} and  \cite{Gur88}.

\smallskip

Implication TT $\Rightarrow$ NS  was shown in \cite{BCG12} for  $n=3$ and
it is an open problem whether it holds for  $n > 3$.
The following three implications hold for all  $n$:

\begin{itemize}
\item DS $\Rightarrow$ T, \cite{Mou83};
\item TT $\Rightarrow$ T, \cite{BGMP10}, by the definitions, it is enough to verify for  $n = 2$;
\item AC $\Rightarrow$ NS  is immediate from the definitions;
\item DS $\Rightarrow$ NS, immediately follows from Proposition \ref{p1}.
\end{itemize}

The example in \cite[Figure 1 in Section 5]{BCG12} disproves
two implications: \;\; TT $\Rightarrow$ AC \;\; and  \;\; TT $\Rightarrow$ DS.
Thus, total tightness becomes much weaker when  $n > 2$.
The following simple example from \cite{Gur88} shows the same for acyclicity.

\bigskip

$a_1 a_1 \;\;\; a_1 a_2$

$a_2 a_1 \;\;\; a_1 a_1$

\bigskip

It is not difficult to verify that for the above game form
AC and, hence, NS  hold, while DS, T, and hence, TT  do not.

\begin{remark}
It was claimed in \cite[Remark 3]{Gur75}
that \;  NS $\Rightarrow$ T.
This mistake was noticed by Danilov.
\end{remark}

Thus, we get an almost complete analysis of
the implications between the above five classes of game forms.
Only \; TT~$\Rightarrow~NS$ \; remains open for   $n > 3$.

\subsection{Dominance-solvable game forms in veto-voting}
\label{ss13}
Here we present briefly the main results of  \cite{GG86}
following the notation of the textbook \cite[Chapter 6]{Mou83}.

In the voting theory, we interpret  $I = \{1, \ldots, n\}$  as voters
and  $A = \{a_1, \ldots, a_p\}$  as candidates.
To each  $i \in I$  we assign a positive integer  $\mu_i$,
interpret it as his veto-power, and give him $\mu_i$   veto-cards.
Respectively, with every  $a \in A$  we associate
a positive integer  $\lambda_a$  that will mean her veto-resistance.
We always assume that  $\sum_{i \in I} \mu_i < \sum_{a \in A}  \lambda_a.$

A strategy  $x_i$  of a voter  $i \in I$  is an arbitrary distribution
of his veto-cards among the candidates  $a \in A$  such that
each one gets at most her veto-resistance, $\mu_i(a) \leq  \lambda_a$.
Given a strategy profile  $x = (x_1, \ldots, x_n)$  we define
$\mu(a) = \sum_{i \in I} \mu_i(a)$.
All candidates  $a \in A$  with  $\mu(a) \geq \lambda_a$  are vetoed,
and all others, $A(x) \subseteq A$  with  $\mu(a) < \lambda_a$,  are elected.
Obviously, at least one candidate is elected, $A(x) \neq \emptyset$, but maybe, more.
Let  $\mu = (\mu_1, \ldots, \mu_n)$  and  $\lambda = (\lambda_1, \ldots , \lambda_p)$
denote veto-power and veto-resistance profiles, respectively.
The mapping  $C = C_{\mu, \lambda} : X \rightarrow 2^A$  defined by  $C(x) = A(x)$
is a veto game correspondence.
Selecting an arbitrary candidate  $a \in A(x)$  for every  $x \in X$
we obtain a veto game form  $g \in C_{\mu, \lambda}$.

\begin{proposition}
A veto game form   $g \in C_{\mu, \lambda}$  is tight whenever
$ \sum_{a \in A}  \lambda_a - \sum_{i \in I} \mu_i = 1$.
\qed
\end{proposition}

The proof can be found in \cite[Chapter 6]{Mou83}.
The result may be easily reformulated in terms
if the so-called threshold Boolean functions \cite[Part II Chapter 9]{CH11}.

\medskip

The DS veto game forms are sparse, however, there are
infinitely many of them \cite{GG86}. Below we provide three examples.

\smallskip
\noindent
{\bf Example 1.}
There are two voters of veto-power $1$  each and
three candidates of veto-resistance  $1$  each, that is,
$n = 2, \mu = (1,1)$ and  $p = 3, \lambda = (1,1,1)$.

\bigskip

$\; a_1 a_2 a_3$  \;\;\;\;\;\;\;\;\;\;\;\;\,  $\; a_1 a_2 a_3$

\smallskip
\noindent
$a_1 \;\;\; a_2 a_3 a_2$ \;\;\;\;\;\;\; $a_1 \;\;\; a_3 a_3 a_2$
\newline
$a_2 \;\;\; a_3 a_3 a_1$ \;\;\;\;\;\;\; $a_2 \;\;\; a_3 a_1 a_1$
\newline
$a_3 \;\;\; a_2 a_1 a_1$ \;\;\;\;\;\;\; $a_3 \;\;\; a_2 a_1 a_2$

\bigskip

Each voter has three strategies: to veto $a_1, a_2$, or $a_3$.
Hence, $C = C_{(1,1), (1,1,1)}$ is a $3 \times 3$ veto game correspondence.
There are $2^3 = 8$  veto game forms  $g \in C$, since
we have to choose one of the two elected candidates on the main diagonal, while
in each of the remaining six strategy profiles the elected candidate is unique.
We obtain a DS  veto game form only when we select
three pairwise distinct candidates on the main diagonal.
There are only two such choices, they are given above.
It is not difficult to verify that both veto game forms are DS.

\begin{remark}
\label{r2}
Let us note that every two strategies
of the voter  $1$  or  $2$  differ by two outcomes
and recall Remark \ref{r1}.
Thus, DS veto game forms are very different from the positional ones.
\end{remark}

\noindent
{\bf Example 2.}
There are two voters of veto-power $1$ and $2$  and
four candidates of veto-resistance  $1$  each, that is,
$n = 2$, $\mu = (1,2)$ and $p = 4, \lambda = (1,1,1,1)$.

\bigskip

$\; a_3 a_2 a_2 a_1 a_1 a_1$  \;\;\;\;\;\;\;\;\;\;\;\;\;\;\;\,       $\; a_3 a_2 a_2 a_1 a_1 a_1$

$\; a_4 a_4 a_3 a_4 a_3 a_2$  \;\;\;\;\;\;\;\;\;\;\;\;\;\;\;\,       $\; a_4 a_4 a_3 a_4 a_3 a_2$

\smallskip
\noindent
$a_4 \;\;\; a_1 a_1 a_4 a_1 a_2 a_3$  \;\;\;\;\;\;\;\;\;\; $a_4 \;\;\; a_1 a_1 a_1 a_2 a_2 a_3$
\newline
$a_3 \;\;\; a_4 a_4 a_4 a_1 a_2 a_3$  \;\;\;\;\;\;\;\;\;\; $a_3 \;\;\; a_1 a_1 a_4 a_2 a_4 a_4$
\newline
$a_2 \;\;\; a_1 a_2 a_3 a_2 a_2 a_3$  \;\;\;\;\;\;\;\;\;\; $a_2 \;\;\; a_1 a_3 a_4 a_3 a_4 a_3$
\newline
$a_1 \;\;\; a_4 a_3 a_3 a_2 a_2 a_3$  \;\;\;\;\;\;\;\;\;\; $a_1 \;\;\; a_2 a_3 a_4 a_2 a_2 a_3$

\bigskip

In this case  $C = C_{(1,2), (1,1,1,1)}$ is a $4 \times 6$ veto game correspondence.
There are $2^{10} = 1024$  veto game forms  $g \in C$ and only four of them are DS \cite{GG86};
both are given above.

\medskip
\noindent
{\bf Example 3.}
There are two voters of veto-power $2$ each and
five candidates of veto-resistance  $1$  each, that is,
$n = 2$, $\mu = (2,2)$ and $p = 4, \lambda = (1,1,1,1,1)$.

\newpage

$\;\;\;\;\, a_4 a_3 a_3 a_2 a_2 a_2 \;\; a_1 a_1 a_1 a_1$

$\;\;\;\;\, a_5 a_5 a_4 a_5 a_4 a_3 \;\; a_5 a_4 a_3 a_2$

\smallskip
\noindent
$a_4 a_5 \;\;\; a_1 a_1 a_1 a_1 a_1 a_1 \;\; a_2 a_2 a_2 a_3$
\newline
$a_3 a_5 \;\;\; a_1 a_1 a_1 a_1 a_1 a_4 \;\; a_2 a_2 a_4 a_4$
\newline
$a_3 a_4 \;\;\; a_1 a_1 a_5 a_1 a_5 a_5 \;\; a_2 a_5 a_5 a_5$
\newline
$a_2 a_5 \;\;\; a_1 a_1 a_1 a_3 a_3 a_4 \;\; a_3 a_3 a_4 a_3$
\newline
$a_2 a_4 \;\;\; a_1 a_1 a_5 a_3 a_5 a_5 \;\; a_3 a_5 a_5 a_3$
\newline
$a_2 a_3 \;\;\; a_1 a_4 a_5 a_4 a_5 a_4 \;\; a_4 a_5 a_4 a_4$

\smallskip
\noindent
$a_1 a_5 \;\;\; a_2 a_2 a_2 a_3 a_3 a_4 \;\; a_2 a_2 a_2 a_3$
\newline
$a_1 a_4 \;\;\; a_2 a_2 a_5 a_3 a_5 a_5 \;\; a_2 a_2 a_2 a_3$
\newline
$a_1 a_3 \;\;\; a_2 a_4 a_5 a_4 a_5 a_4 \;\; a_2 a_2 a_4 a_4$
\newline
$a_1 a_2 \;\;\; a_3 a_4 a_5 a_3 a_3 a_4 \;\; a_3 a_3 a_4 a_3$

\bigskip

In this case  $C = C_{(2,2), (1,1,1,1,1)}$ is a $10 \times 10$ veto game correspondence.
One of its DS veto game forms is given above.

\subsection{On dominance-solvable, sincere, and stable  social choice functions}
\label{ss14}
$\;\;\;\;$ {\em - Hey, you, little gourmand, if you had chosen first,
which nut would you have taken?"

- Don't doubt, I would have taken the small one, -
Junior answered firmly.

- Then why being so worried? You have got it, haven't you?

\smallskip

A. Lindgren. Junior and Karlson who lives on the roof.}

\bigskip

A social choice function (SCF) is a special game form
$g : X \rightarrow A$, where

\bigskip

$A = \{a_1, \ldots, a_p\}, \;  I = \{1, \ldots, n\}; \;
X = \prod_{i \in I} X_i =  X_1 \times \ldots \times X_n$;

\medskip
\noindent
furthermore  $X_1 = \ldots = X_n$  is the set of all linear orders over $A$.
In other words, each voter $i \in I$, as a strategy  $x_i  \in X_i$,
reveals his/her  preference over the set of candidates
$A = \{a_1, \ldots, a_p\}$, and
the elected candidate  $a \in A$  is a function
$a = g(x)$  of the obtained preference profile
$x = (x_1, \ldots, x_n) \in X_1 \times \ldots \times X_n$.
See more details in \cite[Chapter 2]{Mou83}.

An SCF is called {\em sincere}
(or {\em manipulation-free}, or {\em strategy-proof}) if
every strategy profile  $x \in X$   is a NE
with respect to the same preference profile  $x$.
This is an attractive property, since if it holds,
the voters need no manipulations, since
no voter  $i \in I$  can make a profit
revealing a preference  $x'_i$  distinct from his/her actual preference  $x_i$.

Notice that in the above definitions  $x \in X$  is viewed
simultaneously as a strategy and preference profile.

\smallskip

An SCF is called {\em dictatorial} if there is a voter  $i_0 \in I$  such that
in every situation  $x \in X$  the best candidate for  $i_0$  is chosen.

Obviously, the dictator is unique and any such a SCF is sincere.
Indeed, the dictator can only lose by manipulating, while
for any other voter it is just useless,
because his/her preference does not  matter anyway.

Given an  SCF  $S$  with  $n$ voters and only two candidates, $p=2$,
and a profile  $x$,  let   $I = K_1(x) \cap K_2(x)$
denote the partition of the voters into two coalitions
$K_1 = K_1(x)$  and  $K_2 = K_2(x)$, which prefer
candidate  $a_1$  to  $a_2$  and, respectively, vice versa.
Then,  $S$  is called {\em monotone} if for  $j = 1,2$  we have:
$S(x) = a_j \Rightarrow S(x') = a_j$,
provided  $K_j(x) \subseteq K_j(x')$.
In other words, a candidate remain elected
whenever the supporting coalition is getting larger.
Clearly, such SCFs are sincere, too.

\smallskip

By the famous Gibbard-Satterthwaite theorem \cite{Gib73, Sat73},
there are no others. 
This statement is similar and closely related to the
the well-known paradox Arrow in the voting theory;
see \cite[Chapters 4.4 and 3.6]{Mou83} for more details.

\smallskip

Due to this ``negative" result,
it  seems natural to relax the property of sincerity.
We will apply the concept of DS.
If an SCF  $S$  is DS  then a DE  is defined for each  $x \in X$.
Thus, we obtain another SCF  $S' = D(S)$.
Somewhat surprisingly,  $D(S)$  frequently
(although not always) appears in its turn to be DS itself, and
so we can consider the sequence  $S, D(S), D^2(S), \ldots$ etc.,
getting domination cycles.
All such cycles are found in case  $n = p = 3$ \cite{GGM86};
see \cite{GGZ86} for more details.

\medskip

An SCF  $S$  is called {\em stable} if  $D(S) = S$.
A stable  SCF  is not sincere
(unless it is dictatorial or monotone with  $p=2$)
and yet, for every  $x \in X$  the
game  $(S,x)$  in the normal form is DS and, moreover,
its DE outcome  $(D(S))(x) = S'(x)$  is equal to  $S(x)$;
see the epigraph above.

\medskip

Let us recall an example with  $n = p = 3$  from \cite{GGM86}.
Voters $1$ and $2$ are leaders,
together they have full power and may elect
any one of the three candidates.
If they agree upon the best candidate then (s)he is elected,
and  if the worst candidates in their preference lists are distinct
then the remaining one is elected.
Under these two rules, only six conflict situations remain unresolved:
both voters, $1$  and  $2$,  hate most the same candidate
but like most distinct two.
In this case they ask the voter $3$  (who represents people):
``who among the first two candidates do you prefer?"
listen to the answer and choose ... oppositely.
The obtained SCF  $S$  is stable, although, obviously, unfair
(and one can even say ``mocking") with respect to people, represented by voter $3$.
Note that if the two leaders would choose
in accordance with the people's will then
the obtained SCF  $S'$  won't be stable; see more details in  \cite{GGM86, GGZ86}.

\subsection{Criteria of Nash- and dominance-solvability
based on excluding $2 \times 2$-subconfigurations}
\label{ss15}
\subsubsection{Recognition of math-patters}
\label{sss150}
In 1964 Shapley \cite{Sha64}  observed that a matrix 
has a saddle point (SP) whenever every its $2 \times 2$  submatrix has one.
In other words, the real $2 \times 2$ matrices in which
both entries of one diagonal
is strictly larger than both entries of the other
are the only minimal SP-free matrices.
This result was strengthen in \cite{BGM09}, where it was demonstrated
that these matrices are not just the only minimal,
but also the only locally minimal SP-free matrices.
In other words, every SP-free matrix of size larger than $2 \times 2$
has a row or a column that can be deleted and
the reduced matrix still remains SP-free.

This property can be extended to the two-person non-zero-sum case, but
in a more sophisticated way.
A locally minimal NE-free bimatrix games were
fully characterized in \cite{BGM09};
see also \cite[Theorem 1 in Section 1.2]{BEGMO16}.
Such a bimatrix is square, but may be of any size  $n \times n$,  with  $n \geq 2$,
and is representable as the direct sum of strong improvement cycles;
see \cite{BGM09} for the definitions and more details.

Furthermore, Shapley's \cite{Sha64} statement can be generalized
to bimatrix games in many other ways as follows.
Recall that, in general, NE, DE, and the concept of domination itself
depend on the pre-orders of the players rather than on the real values
of the utility function.
So, let us partition all $2 \times 2$ bimatrix games into fifteen classes
$S = \{c_1, \ldots, c_{15}\}$  depending on the preference pre-orders of the two players.

A subset  $t \subseteq S$  is called a DE-theorem
if a bimatrix game has a DE,
whenever it contains no $2 \times 2$ subgame from  $t$.
A subset  $e \subseteq S$  is called a DE-example
if a bimatrix game has no DE, while
all its  $2 \times 2$ subgames are from  $e$.
It is not difficult to show that two set-families
$T_{DE}$ and  $E_{DE}$  of all minimal
(that is, the strongest)
DE-theorems and DE-examples are dual
(sometimes called also transversal).
This property provides an efficient stopping criterion:
we extend the considered two lists verifying each time the duality condition
which will certify that both lists are complete and cannot be extended further.
By this method we obtain all DE-theorems and all DE-examples.
One can replace DE by NE, the method remains the same.

Actually, this method, which
can be called {\em the recognition of math-patterns}, is even much more general.
One can consider any set of attributes (properties) $S$
trying to characterize any target subset $P_0$
within an arbitrary set of objects  $P$.
In \cite{GG83}, this approach was illustrated
by a simple model problem in which
$P$  is the set of all $4$-gons on the plane,
$P_0$  is the subset of squares, and
$S$  is a set of eight properties of a quadrilateral. 
Two dual hypergraphs of all minimal theorems
$T_{SQ}$  and examples  $E_{SQ}$  were constructed.
In \cite{CLS90}, this approach was applied
to a more serious problem related to
some special families of the Berge graphs.
Recently, the method was analyzed in \cite{BEGMO16}
to several concepts of game theory.
Here we will recall two cases, related to NE and DE for bimatrix games.
So from now on we restrict ourselves to the case of two players, $n=2$.

\begin{figure}[t]
\centerline{
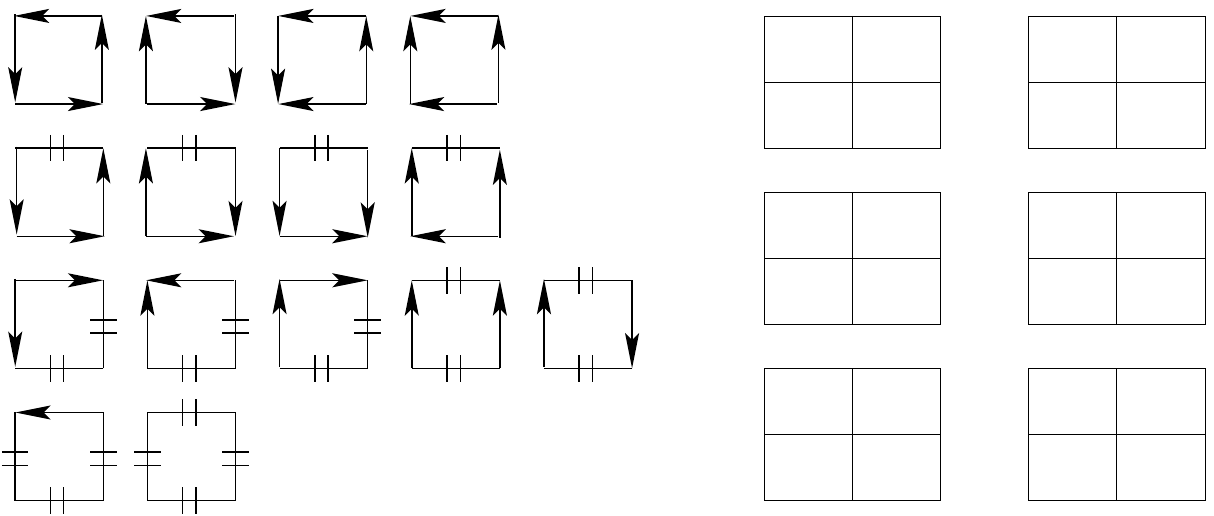
}
\caption{Fifteen $2$-squares.
\newline
In Figures 1--6  we use the following notation.
Configurations are represented by planar grids whose
nodes correspond to the strategy profiles  $(x_1, x_2) \in X_1 \times X_2$.
A line with two dashes between two nodes means
that the corresponding two situations make a tie, while
an arrow from one node to another means that
the second situation is better than the first one,
for the player choosing among them.
The corresponding games are represented by tables whose rows and columns are
the strategies  $X_1$  and  $X_2$, respectively.
Thus, strategy profiles  $(x_1, x_2) \in X_1 \times X_2$  are represented by
the cells of these tables, where
$u_1(x_1, x_2)$  and  $u_2(x_1, x_2)$  are located respectively
in the bottom-left and top-right corners of the cell  $(x_1, x_2)$.
}
\label{f1}
\end{figure}

\subsubsection{Configurations, alphabet, fifteen $2$-squares}
\label{sss151}
To decide whether a strategy profile
$x = (x_1, x_2) \in X_1 \times X_2 = X$  is a NE
or whether  $x_i \succeq x'_i$  for  $i \in \{1,2\}$,
One does not need to know the payoff function $u : X \rightarrow R^2$.
It is enough to know pre-orders:
for any  $x_i, x'_i \in X_i$  and  $x_{3-i} \in X_{3-i}$   whether
$i$  prefers  $(x_i, x_{3-i})$  to   $(x'_i, x_{3-i})$,
or vice versa, or (s)he is indifferent and these two
strategy profiles make a tie.

\smallskip 

A product  $X = X_1 \times X_2$  with all
preference pre-orders of both players
will be called a {\em configuration};
restricted to a box  $X'_1 \times X'_2$, where
$X'_i \subseteq X_i$  for  $i =1,2$, it
will be called a {\em subconfiguration}.
For brevity, a subconfiguration of size  $2 \times 2$,
that is, when  $|X_1| = |X_2| = 2$, is called a {\em 2-square}.
We will not distinguish $2$-squares one of which can be obtained
from the other by a permutation of strategies or by the transposition.
Then, it is easy to verify that there exist exactly fifteen $2$-squares
given in Figure \ref{f1}.
Among these fifteen $2$-squares, fourteen have an NE, all except  $c_1$, and
thirteen have a DE, all except  $c_1$ and $c_2$.

Note that the first six $2$-squares are frequent in the textbooks.
For example, $c_2$  and  $c_4$  represent the classical
``battle of sexes"  and ``prisoner's dilemma",
while  $c_5$  and  $c_6$  illustrate the concepts
of the promise and threat, respectively.
The corresponding bimatrix games are given in Figure 1.

Shapley's theorem asserts that every  $c_1$-free zero-sum game
(or configuration) has a saddle point.
Thus, it outlines a hereditary class of NS zero-sum games. 
In this section we obtain similar conditions
for NS and DS, in terms of some subsets of forbidden $2$-squares
from $S = \{c_1, \ldots, c_{15}\}$. 
Obviously, any class of games
(more generally, objects) defined by a family of forbidden subgames
(respectively, subobjects) is hereditary.  
In contrast, the families of the (zero-sum) NS or DS games are not.   
Hence, Shapley's and all other conditions announced above 
are only sufficient but not necessary.
They define some hereditary subclasses of NS and DS  games.

\subsubsection{NE-examples, NE-conjectures, and NE-theorems}
\label{sss152}
Our goal is to find and prove all strongest NE-theorems, which are assigned to
all inclusion-minimal subsets  $t \subseteq S = \{c_1, \ldots, c_{15}\}$
such that any $t$-free configuration contains an NE.
To this end, we start to collect the strongest NE-examples,
that is, all inclusion-minimal subsets  $e \subseteq S$ such that
there exists an NE-free configuration containing only $2$-squares from  $e$.
The simplest such configuration is the $2$-square  $c_1$  alone; $e = (c_1)$.
Four others are shown in Figure 2.
Thus, we obtain a family of NE-examples
$\{(c_1), (c_2, c_3), (c_3, c_5, c_6), (c_5, c_9), (c_2, c_4, c_5, c_6)\}$.
It is convenient to write it as a disjunctive normal form (DNF)

\medskip

$E_{NE} = c_1 \vee c_2 c_3 \vee c_2 c_4 c_5 c_6 \vee c_3 c_5 c_6 \vee c_5 c_9$.

\begin{remark}
As we already mentioned in the beginning of Section \ref{ss15},
a locally minimal NE-free bimatrix games
(and configurations, as well) were fully characterized in \cite{BGM09}.
They may be of any size  $n \times n$,  with  $n \geq 2$,
and representable as direct sums of strong improvement cycles;
see \cite{BGM09} for the definition more details.
\newline
This result simplifies constructing NE-examples a lot.
Obviously, each NE-example is a locally minimal NE-free configuration
and, hence, it must contain a strong improvement cycles;
$c_1$  is such a cycle itself.
Each configuration in Figure 2
is a square of size $n \times n$,  for some  $n \geq 2$,
containing a strong improvement cycle  $C_n$.
The corresponding  $2n$  entries are denoted by
$n$  white and  $n$  black small disks such that
a white disk is a unique maximum in its row for the column player, while
a black disk is a unique maximum in its column for the row player.
\end{remark}

\begin{figure}[t]
\centerline{
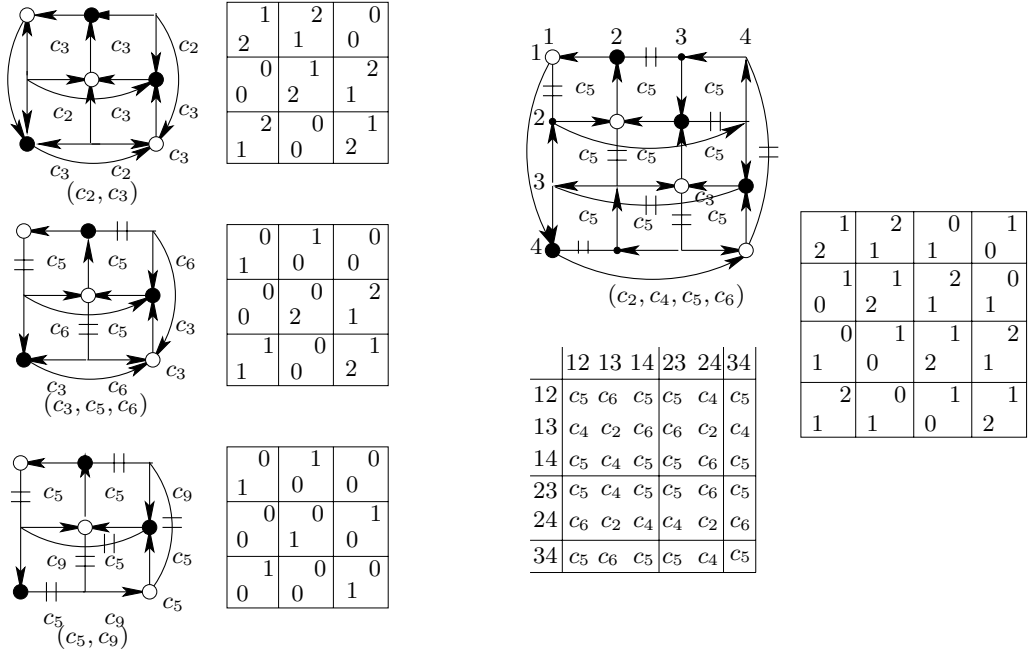
}
\caption{NE-examples.
Four NE-free games and the corresponding configurations
whose sets of types of $2$-cycles,
$\{(c_2, c_3), (c_3, c_5, c_6), (c_5, c_9)$, and $(c_2, c_4, c_5, c_6)\}$,
form all NE-examples.
\newline
In Figures \ref{f2}--6 we keep the same legend as for Figure \ref{f1}.
Yet, two nodes may be not connected when
the relation between the corresponding two strategy profiles follows from transitivity.
}
\label{f2}
\end{figure}

It is clear that every NE-theorem $t \in T_{NE}$  and
NE-example  $e \in E_{NE}$  must intersect.
Indeed, if $e \cap t \neq \emptyset$  then
$e$  is a counterexample to  $t$  and, hence, it is not a theorem.
Moreover, our research is complete whenever $E_{NE}$  and $T_{NE}$
form a pair of dual hypergraphs, since in this case,
by the definition of duality, any new example or theorem would be weaker
than one that already exists.
(This is true not only for NE, but for any target property
$P_0 \subseteq P$  and for any alphabet $S$ of arbitrary attributes.)
Dualizing DNF  $E_{NE}$  we obtain

$$E^d_{NE} = c_1 (c_2 c_5 \vee c_3 c_5  \vee c_2 c_3 c_9 \vee c_3 c_4 c_9  \vee c_3 c_6 c_9) =
 c_1 c_2 c_5 \vee c_1 c_3 c_5  \vee c_1 c_2 c_3 c_9 \vee c_1 c_2 c_6 c_9 \vee c_1 c_3 c_4 c_9  \vee c_1 c_3 c_6 c_9.$$

In \cite{BEGMO16} all six theorems were proven.
Thus, $T_{NE} = E^d_{NE}$  and our NE-research is complete,
in terms of the set of properties $S = \{c_1, \dots, c_{15}\}$.
Of course,  $S$ can  be modified.
For example, it is not mandatory to identify the transposed $2$-squares,
we can view them as distinct.
Doing so, we expand the set of attributes and change
the families of NE-examples and NE-theorems; although they will remain dual.

Let us underline that some examples may be missing in $E$
and, then, some theorems may be missing in $E^d$;
moreover, some ``theorems" from  $E^d$  may be just wrong.
Suppose that we failed to construct the example
$e = (c_2, c_4, c_5, c_6) \in E_{NE}$,  which
is possible, since the corresponding $4 \times 4$  configuration is not that simple.
Then dualizing the ``incomplete" DNF
$E'_{NE} = c_1 \vee c_2 c_3 \vee c_3 c_5 c_6 \vee c_5 c_9$, we obtain

$$E^d_{NE} = c_1 (c_2 c_5 \vee c_3 c_5  \vee c_2 c_6 c_9 \vee c_3 c_9) =
 c_1 c_2 c_5 \vee c_1 c_3 c_5  \vee c_1 c_2 c_6 c_9 \vee c_1 c_3 c_9.$$

Instead of three ugly theorems
$c_1 c_2 c_3 c_9 \vee c_1 c_3 c_4 c_9  \vee c_1 c_3 c_6 c_9$, we have
one nice conjecture $c_1 c_3 c_9$, which has only one disadvantage:
it is an overstatement.
If we fail proving it after sufficiently hard work and long time,
we will have to return to our list of examples and verify it.
Perhaps, we can strengthen some of them or add one
(as in the considered case) or several missing examples.
If we cannot do this, we return to the list of theorems again, etc.
The result is not guaranteed, as subsection \ref{sss154} shows.

\subsubsection{DE-examples, DE-conjectures, and DE-theorems}
\label{sss153}
Let us note that  $c_1$  must be a DE-example,
because it is a NE-free $2$-square.
Furthermore, $c_2$  is a DE-example too, although it has a NE.
Except these two, we got five larger DE-examples:
$(c_3, c_5, c_6)$ from Figure 2 and four shown in Figures 3--6.

\begin{remark}
Note that  $(c_2, c_4, c_5, c_6)$  is
a DE-example too, since it is an NE-example, see Figure 2,
but it is not minimal, since  $(c_4, c_5, c_6)$
is a stronger DE-example; see Figure 6. Yet, the latter has a NE.
\end{remark}

Thus, we obtain the DNF
$E_{DE} = c_1 \vee c_2 \vee c_3 c_5 c_6 \vee c_4 c_5 c_6 \vee c_5 c_9 \vee c_5 c_{10} \vee  c_6 c_{11}.$
Dualizing it, we  get the list of DE-conjectures,
which all can be proven.
Thus, $E^d_{DE} = T_{DE}$  and we obtain the complete DE-research:

$$E^d_{DE} = c_1 c_2 (c_3 c_4 c_9 c_{10} c_{11} \vee c_5 c_6 \vee c_5 c_{11} \vee c_6 c_9 c_{10}) =
c_1 c_2 c_3 c_4 c_9 c_{10} c_{11} \vee c_1 c_2 c_5 c_6 \vee c_1 c_2 c_5 c_{11} \vee c_1 c_2 c_6 c_9 c_{10}.$$

Note that the last two DE-theorems,
$(c_1, c_2, c_5, c_{11})$  and  $(c_1, c_2, c_6, c_9, c_{10})$,
are equivalent, since one is obtained from the other
by changing the pre-order of a player
(one from two) to the inverse one.
Obviously, this operation keeps all domination free subgames and
transforms one of the above two DE-theorems into the other.

Even more important is to notice that
in this subsection we waved the assumption
$u_1(x) = u_1(x') \Leftrightarrow u_2(x) = u_2(x')$
and, as we know, without this assumption
the concept of domination is not well-defined, although
formally all above results are correct.
For this reason, in the next subsection
we modify the NE- and DE-research adapting
it for the game forms rather than the games in normal form.

\begin{figure}
	\centering
	\includegraphics[width=5in]{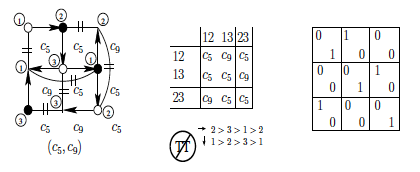}
	\caption{A domination-free non-TT configuration.}
\end{figure}

\begin{figure}
	\centering
	\includegraphics[width=5in]{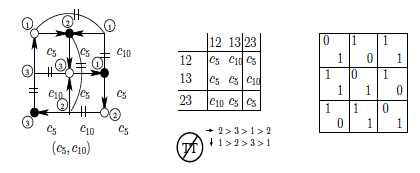}
	\caption{A domination-free non-TT configuration.}
\end{figure}

\begin{figure}
	\centering
	\includegraphics[width=5in]{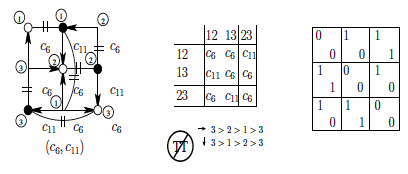}
	\caption{A domination-free non-TT configuration.}
\end{figure}


\begin{figure}
	\centering
	\includegraphics[width=6in]{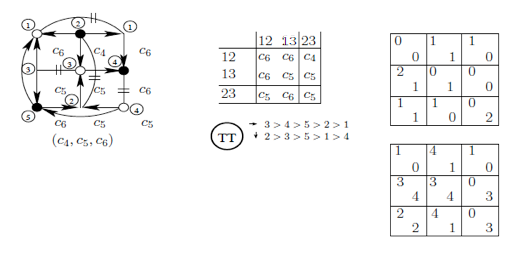}
	\caption{A domination-free TT configuration; the entry $(2,3)$ is an NE}
\end{figure}

\subsubsection{Tie-transitive (or formal) Nash- and dominance-solvability}
\label{sss154}
A configuration is called {\em tie-transitive} (TT) or {\em formal}
if it can be realized by a game form.
One can easily verify in linear time
whether a given configuration is TT.
The transitive closure of all ties
(equalities) of both players uniquely defines a game form.
Extend the pre-order of the player  $i = 1,2$
to the outcomes of this game form.
The configuration is TT  if and only if
no cycle appears, neither for player $1$  nor for $2$.

For example, among fifteen $2$-squares on Figure 1, all are TT,
except for  $c_{13}$  and  $c_{14}$.
It is also not difficult to verify that
the first two configurations on Figure 2 are TT,
while the last two are not.
Thus, in the TT case we keep NE-examples
$E^{TT}_{NE} =  c_1 \vee c_2 c_3 \vee c_3 c_5 c_6$,
which still can be found in Figure 2.
Dualizing we obtain the NE-conjectures
$E^{TT \;\; d}_{NE} =  c_1 c_2 c_5 \vee c_1 c_3 c_5  \vee c_1 c_3$,
all of which appear to be NE-theorems.
Thus, we obtain the perfect NE-research.

\medskip

However, problems appear with domination.
From Figures 2 and 6 we obtain
$E^{TT}_{DE} =  c_1 \vee c_2  \vee c_3 c_5 c_ 6 \vee c_4 c_5 c_6$.
Dualizing we obtain only three DE-conjectures

$$E^{TT \;\; d}_{DE} =  c_1 c_2 (c_3 c_4 \vee c_5 \vee c_6) =
 c_1 c_2 c_3 c_4 \vee c_1 c_2 c_5 \vee c_1 c_2 c_6,$$
\noindent
from which only the first one is proven, while the last
two remain (equivalent) conjectures.
It is also possible that some minimal TT DE-example is missing.

\begin{remark}
The TT DE-example  $(c_4, c_5, c_6)$  on Figure 6 is ``stronger" than
minimal NE-example $(c_2, c_4, c_5, c_6)$ in Figure 2.
Moreover, the first one is TT, while the second one is not.
However, the first one has a NE, although it is domination-free.
\end{remark}


\end{document}